\documentclass[a4paper,twoside,10pt]{amsart}
\usepackage{amsmath,amsfonts,amssymb,amsthm}
\usepackage{mathrsfs}
\usepackage{enumerate}
\newtheorem{theorem}{Theorem}[section]
\newtheorem{corollary}{Corollary}[section]
\numberwithin{equation}{section}

\author[G. Nemes]{Gerg\H{o} Nemes}
\address{Central European University, Department of Mathematics and its Applications, H-1051 Budapest, N\'ador utca 9, Hungary}
\email{nemesgery@gmail.com}

\begin{document}

\title[An explicit formula for coefficients]{An explicit formula for\\ the coefficients in Laplace's method}

\begin{abstract}
Laplace's method is one of the fundamental techniques in the asymptotic approximation of integrals. The coefficients appearing in the resulting asymptotic expansion, arise as the coefficients of a convergent or asymptotic series of a function defined in an implicit form. Due to the tedious computation of these coefficients, most standard textbooks on asymptotic approximations of integrals do not give explicit formulas for them. Nevertheless, we can find some more or less explicit representations for the coefficients in the literature: Perron's formula gives them in terms of derivatives of an explicit function; Campbell, Fr\"oman and Walles simplified Perron's method by computing these derivatives using an explicit recurrence relation. The most recent contribution is due to Wojdylo, who rediscovered the Campbell, Fr\"oman and Walles formula and rewrote it in terms of partial ordinary Bell polynomials. In this paper, we provide an alternative representation for the coefficients, which contains ordinary potential polynomials. The proof is based on Perron's formula and a theorem of Comtet. The asymptotic expansions of the gamma function and the incomplete gamma function are given as illustrations.
\end{abstract}
\maketitle

\section{Introduction}

Laplace's method is one of the best-known techniques developing asymptotic approximation for integrals. The origins of the method date back to Pierre-Simon de Laplace (1749 -- 1827), who studied the estimation of integrals, arising in probability theory, of the form
\begin{equation}\label{eq1}
I\left( \lambda  \right) = \int_a^b {\mathrm{e}^{ - \lambda f\left( x \right)} g \left( x \right)\mathrm{d}x} \quad \left(\lambda \to +\infty\right).
\end{equation}
Here $\left(a,b\right)$ is a real (finite or infinite) interval, $\lambda$ is a large positive parameter and the functions $f$ and $g$ are continuous. Laplace made the observation that the major contribution to the integral $I\left( \lambda  \right)$ should come from the neighborhood of the point where $f$ attains its smallest value. Observe that by subdividing the range of integration at the minima and maxima of $f$, and by reversing the sign of $x$ whenever necessary, we may assume, without loss of generality, that $f$ has only one minimum in $\left[a, b\right]$ occuring at $x = a$. With certain assumptions on $f$, Laplace's result is
\[
\int_a^b {\mathrm{e}^{ - \lambda f\left( x \right)} g \left( x \right)\mathrm{d}x} \sim g\left( {a} \right)\mathrm{e}^{ - \lambda f\left( {a} \right)} \sqrt {\frac{{\pi }}{{2 \lambda f''\left( {a} \right)}}} .
\]
The sign $\sim$ is used to mean that the quotient of the left-hand side by the right-hand side approaches $1$ as $\lambda \to +\infty$.  This formula is now known as Laplace's approximation. A heuristic proof of this formula may proceed as follows. First, we replace $f$ and $g$ by the leading terms in their Taylor series expansions around $x = a$, and then we extend the integration limit to $+\infty$. Hence,
\begin{align*}
\int_a^b {\mathrm{e}^{ - \lambda f\left( x \right)} g\left( x \right)\mathrm{d}x} & \approx \int_a^b {\mathrm{e}^{ - \lambda \left( {f\left( {a} \right) + \frac{{f''\left( {a} \right)}}{2}\left( {x - a} \right)^2 } \right)} g\left( {a} \right)\mathrm{d}x} \\
& \approx g\left( {a} \right)\mathrm{e}^{ - \lambda f\left( {a} \right)} \int_{a}^{ + \infty } {\mathrm{e}^{ - \lambda \frac{{f''\left( {a} \right)}}{2}\left( {x - a} \right)^2 } \mathrm{d}x} \\
& = g\left( {a} \right)\mathrm{e}^{ - \lambda f\left( {a} \right)} \sqrt {\frac{{\pi }}{{2 \lambda f''\left( {a} \right)}}} .
\end{align*}
The modern version of the method was formulated in 1956 by Arthur Erd\'elyi (1908 -- 1977), who applied Watson's lemma to obtain a complete asymptotic expansion for the integral \eqref{eq1}. His method requires some assumptions on $f$ and $g$. Suppose, again, that $f$ has only one minimum in $\left[a, b\right]$ which occurs at $x = a$. Assume also that, as $x \to a^+$,
\begin{equation}\label{exp1}
f\left( x \right) \sim f\left( a \right) + \sum\limits_{k = 0}^\infty {a_k \left( {x - a} \right)^{k + \alpha } } ,
\end{equation}
and
\begin{equation}\label{exp2}
g \left( x \right) \sim \sum\limits_{k = 0}^\infty {b_k \left( {x - a} \right)^{k + \beta  - 1} }
\end{equation}
with $\alpha>0$, $\Re \left( \beta  \right) > 0$, and that the expansion of $f$ can be term-wise differentiated, that is,
\begin{equation}\label{exp3}
f'\left( x \right) \sim \sum\limits_{k = 0}^\infty {a_k \left( {k + \alpha } \right)\left( {x - a} \right)^{k + \alpha  - 1} } 
\end{equation}
as $x \to a^+$. We may also assume, without loss of generality, that $a_0  \ne 0$ and $b_0  \ne 0$. The following theorem is Erd\'elyi's formulation of Laplace's classical method.

\begin{theorem}\label{theorem1} For the integral
\[
I\left( \lambda  \right) = \int_a^b {\mathrm{e}^{ - \lambda f\left( x \right)} g \left( x \right)\mathrm{d}x} ,
\]
we assume that
\begin{enumerate}[(i)]
	\item $f(x)>f(a)$ for all $x \in (a,b)$, and for every $\delta>0$ the infimum of $f(x)-f(a)$ in $\left[a+\delta,b\right)$ is positive;
	\item $f'(x)$ and $g(x)$ are continuous in a neighborhood of $x=a$, expect possibly at $a$;
	\item the expansions \eqref{exp1}, \eqref{exp2} and \eqref{exp3} hold; and
	\item the integral $I\left( \lambda  \right)$ converges absolutely for all sufficiently large $\lambda$.
\end{enumerate}
Then
\begin{equation}\label{eq2}
I\left( \lambda  \right) \sim \mathrm{e}^{ - \lambda f\left( a \right)} \sum\limits_{n = 0}^\infty {\varGamma \left( {\frac{{n + \beta }}{\alpha }} \right)\frac{{c_n }}{{\lambda ^{\left( {n + \beta } \right)/\alpha } }}} ,
\end{equation}
as $\lambda \to +\infty$, where the coefficients $c_n$ are expressible in terms of $a_n$ and $b_n$.
\end{theorem}

The first three coefficients $c_n$ are given explicitly by
\[
c_0  = \frac{{b_0 }}{{a_0^{\beta /\alpha } \alpha }},\quad \quad c_1  = \frac{1}{{a_0^{\left( {\beta  + 1} \right)/\alpha } }}\left( {\frac{{b_1 }}{\alpha } - \frac{{\left( {\beta  + 1} \right)a_1 b_0 }}{{\alpha ^2 a_0 }}} \right),
\]
and
\[
c_2  = \frac{1}{{a_0^{\left( {\beta  + 2} \right)/\alpha}}}\left( {\frac{{b_2 }}{\alpha } - \frac{{\left( {\beta  + 2} \right)a_1 b_1 }}{{\alpha ^2 a_0 }} + \left( {\left( {\beta  + \alpha  + 2} \right)a_1^2  - 2\alpha a_0 a_2 } \right)\frac{{\left( {\beta  + 2} \right)b_0 }}{{2\alpha ^3 a_0^2 }}} \right).
\]
To understand the origin of these coefficients, we sketch the proof of Erd\'elyi's theorem. Detailed proofs are given in many standard textbooks on asymptotic analysis, e.g., in Erd\'elyi's original monograph \cite[p. 38]{Erdelyi} or the classical books of Olver \cite[p. 81]{Olver2} and Wong \cite[p. 58]{Wong}. By conditions (ii) and (iii) there exists a number $c \in (a,b)$ such that $f'(x)$ and $g(x)$ are continuous in $(a,c]$, and $f'(x)$ is also positive in this interval. Let $T=f(c)-f(a)$, and define the new variable $t=t(x)$ as
\[
t=f(x)-f(a).
\]
Since $f(x)$ is increasing in $(a,c)$, we can write
\begin{equation}\label{eq3}
\int_a^c {\mathrm{e}^{ - \lambda f\left( x \right)} g\left( x \right)\mathrm{d}x}  = \mathrm{e}^{-\lambda f\left( a \right)} \int_0^T {\mathrm{e}^{ - \lambda t} h\left( t \right)\mathrm{d}t} 
\end{equation}
with $h(t)$ being the continuous function in $\left(0,T\right]$ given by
\begin{equation}\label{eq4}
h(t) = g(x)\frac{\mathrm{d}x}{\mathrm{d}t}= \frac{g(x)}{f'(x)}.
\end{equation}
By assumption,
\[
t \sim \sum\limits_{k = 0}^\infty {a_k \left( {x - a} \right)^{k + \alpha } } \quad\text{as } x\to a^+.
\]
And so, by series reversion, we obtain
\[
x - a \sim \sum\limits_{k = 1}^\infty {d_k t^{k/\alpha } }\quad\text{as } t\to 0^+ .
\]
Substituting this into \eqref{eq4} and using the asymptotic expansions \eqref{exp2}-\eqref{exp3} yields
\[
h\left( t \right) \sim \sum\limits_{k = 0}^\infty {c_k t^{\left( {k + \beta } \right)/\alpha  - 1} } 
\]
as $t\to 0^+$. (In the common case where $g \equiv 1$, we have $\beta = 1$ and $c_k = \left(k + 1\right)d_{k+1}/\alpha$.) We now apply Watson's lemma to the integral on the right-hand side of \eqref{eq3} to obtain
\[
\int_a^c {\mathrm{e}^{ - \lambda f\left( x \right)} g\left( x \right)\mathrm{d}x} \sim \mathrm{e}^{ - \lambda f\left( a \right)} \sum\limits_{n = 0}^\infty {\varGamma \left( {\frac{{n + \beta }}{\alpha }} \right)\frac{{c_n }}{{\lambda ^{\left( {n + \beta } \right)/\alpha } }}} ,
\]
as $\lambda \to +\infty$. To complete the proof, one shows that the integral on the remaining range $(c,b)$ is negligible.

The coefficients $c_n$ are traditionally calculated manually in the way sketched above: by computing the reversion coefficients $d_n$, substituting the resulting series into the asymptotic expansion of \eqref{eq4}, expanding the series in powers of $t$, and equating the coefficients in each specific application of Laplace's method. Due to this somewhat tedious computation, most standard textbooks on asymptotic approximations of integrals do not give explicit formulas for the $c_n$'s. Nevertheless, there are certain formulas in different degrees of explicitness for the coefficients $c_n$ in the literature. When $\left( {x - a} \right)^{ - \alpha } \left( {f\left( x \right) - f\left( a \right)} \right)$ and $\left( {x - a} \right)^{1 - \beta } g\left( x \right)$ are analytic at $x=a$, Perron's formula gives the coefficients in terms of derivatives of an explicit function involving $f$ and $g$. Campbell, Fr\"oman and Walles rediscovered Perron's method and went further by computing these derivatives using an explicit recurrence formula \cite{Campbell}. The most recent contribution is given by Wojdylo, who rediscovered the Campbell, Fr\"oman and Walles formula and rewrote it in terms of partial ordinary Bell polynomials \cite{Wojdylo1}\cite{Wojdylo2}. Using new ideas of combinatorial analysis, he was able to simplify and systematize the computation of the $c_n$'s. The definition of the partial ordinary Bell polynomials provides a direct way of expanding the higher derivatives in Perron's formula. However, there is a formula by Comtet allowing us to obtain a representation for those derivatives in terms of ordinary potential polynomials. This gives a new and alternative method for the computation of the coefficients $c_n$ (see Corollary \ref{maintheorem}). We discuss these formulas in details in Section \ref{section2}.

In Section \ref{section3}, we give two illustrative examples to demonstrate the application of our new method and to compare the results with those given by Wojdylo's formula. In the first example, we obtain several explicit expressions for the Stirling coefficients appearing in the asymptotic expansion of the gamma function. In the second example, we investigate certain polynomials related to the coefficients in the uniform asymptotic expansion of the incomplete gamma function.

The definition and basic properties of the partial ordinary Bell polynomials and the ordinary potential polynomials are collected in Appendix \ref{appendix}. For further discussion on these polynomials, we refer the reader to the book of Comtet \cite[pp. 133--153]{Comtet} or Riordan \cite[pp. 189--191]{Riordan}.

\section{Explicit formulas for the coefficients $c_n$}\label{section2}

The asymptotic theory of integrals of type \eqref{eq1} is also well established when $\lambda$ is complex and $f$, $g$ are holomorphic functions in a domain of the complex plane containing the path of integration $\mathscr{C}$ joining $a$ to $b$. A well-known method for obtaining asymptotic expansions for such integrals is the method of steepest descents (for a detailed discussion of this method, see, e.g., \cite[pp. 5--99]{Paris2} or \cite[pp. 84--103]{Wong}). This method requires the deformation of $\mathscr{C}$ into a specific path that passes through one or more saddle points of $f$ such that the function $\Im\left(f\right)$ is constant on it. (Recall that $z_0$ is a saddle point of $f$ iff $f'\left(z_0\right)=0$.) This new path is called the path of steepest descent. However, in many specific cases, the construction of such a path can be extremely complicated. This problem may be bypassed using Perron's method which -- by requiring some extra assumptions -- avoids the computation of the path of steepest descent, and provides an explicit expression for the coefficients in the resulting asymptotic series \cite{Perron}\cite[p. 103]{Wong}. A direct adaptation of Erd\'{e}lyi's theorem to complex integrals was formulated by Olver, where the coefficients in the asymptotic expansion can be computed formally in the same way as in the real case \cite[p. 121]{Olver2}. In addition, both Perron's and Olver's method allow the functions $f$ and $g$ to have algebraic singularities at the endpoint $a$ with convergent expansions of the form \eqref{exp1} and \eqref{exp2}.

Based on Wojdylo's work, we show below that Perron's explicit formula also holds formally in the case of Erd\'{e}lyi's theorem, when the expansions \eqref{exp1} and \eqref{exp2} may be merely asymptotic. Starting from the equation $h\left(t\right)\mathrm{d}t = g\left(x\right)\mathrm{d} x$, Wojdylo showed that
\[
c_n^\ast \stackrel{\mathrm{def}}{=} \frac{{a_0^{\left( {n + \beta } \right)/\alpha } \alpha }}{{b_0 }}c_n  = \sum\limits_{k = 0}^n {\frac{{b_{n - k} }}{{b_0 }}\frac{1}{{k!}}\left[ {\frac{\mathrm{d}^k}{\mathrm{d}x^k}\left( {1 + \sum\limits_{k = 1}^\infty  {\frac{{a_k }}{{a_0 }}x^k } } \right)^{ - \left( {n + \beta } \right)/\alpha } } \right]_{x = 0} },
\]
where $\mathrm{d}^k / \mathrm{d}x^k$ is the formal $k$th derivative with respect to $x$. Through his analysis, he used these scaled coefficients $c_n^\ast$. We shall not, however, use them in our paper. His expression can be simplified to
\begin{align}
c_n & = \frac{1}{{\alpha a_0^{\left( {n + \beta } \right)/\alpha } }}\sum\limits_{k = 0}^n {\frac{{b_{n - k} }}{{k!}}\left[ {\frac{{\mathrm{d}^k }}{{\mathrm{d}x^k }}\left( {1 + \sum\limits_{k = 1}^\infty  {\frac{{a_k }}{{a_0 }}x^k } } \right)^{ - \left( {n + \beta } \right)/\alpha } } \right]_{x = 0} } \label{eq5} \\
& = \frac{1}{{\alpha a_0^{\left( {n + \beta } \right)/\alpha } }}\sum\limits_{k = 0}^n {\frac{{b_{n - k} }}{{k!}}\left[ {\frac{{\mathrm{d}^k }}{{\mathrm{d}x^k }}\left( {\frac{{a_0 x^\alpha  }}{{\sum\nolimits_{k = 0}^\infty  {a_k x^{k + \alpha } } }}} \right)^{\left( {n + \beta } \right)/\alpha } } \right]_{x = 0} }\nonumber \\
& = \frac{1}{{\alpha a_0^{\left( {n + \beta } \right)/\alpha } }}\sum\limits_{k = 0}^n {\frac{{b_{n - k} }}{{k!}}\left[ {\frac{{\mathrm{d}^k }}{{\mathrm{d}x^k }}\left( {\frac{{a_0 \left( {x - a} \right)^\alpha  }}{{\sum\nolimits_{k = 0}^\infty  {a_k \left( {x - a} \right)^{k + \alpha } } }}} \right)^{\left( {n + \beta } \right)/\alpha } } \right]_{x = a} } .\nonumber
\end{align}
If we identify $f$ with its asymptotic expansion \eqref{exp1}, we arrive at the final forms
\begin{gather}\label{eq6}
\begin{split}
c_n & = \frac{1}{{\alpha a_0^{\left( {n + \beta } \right)/\alpha } }}\sum\limits_{k = 0}^n {\frac{{b_{n - k} }}{{k!}}\left[ {\frac{{\mathrm{d}^k }}{{\mathrm{d}x^k }}\left( {\frac{{a_0 \left( {x - a} \right)^\alpha  }}{{f\left( x \right) - f\left( a \right)}}} \right)^{\left( {n + \beta } \right)/\alpha } } \right]_{x = a} } \\
& = \frac{1}{{\alpha n!a_0^{\left( {n + \beta } \right)/\alpha } }}\left[ {\frac{{\mathrm{d}^n }}{{\mathrm{d}x^n }}\left\{ {G\left( x \right)\left( {\frac{{a_0 \left( {x - a} \right)^\alpha  }}{{f\left( x \right) - f\left( a \right)}}} \right)^{\left( {n + \beta } \right)/\alpha } } \right\}} \right]_{x = a} ,
\end{split}
\end{gather}
where $G\left(x\right)$ is the formal power series $\sum\nolimits_{k = 0}^\infty  {b_k \left( {x - a} \right)^k }$. An alternative form of this expression -- using the formal residue operator -- is stated in \cite[p. 45]{NIST}. If $\left( {x - a} \right)^{ - \alpha } \left( {f\left( x \right) - f\left( a \right)} \right)$ and $\left( {x - a} \right)^{1 - \beta } g\left( x \right)$ are analytic at $x=a$, this is the Perron formula. Applying the definition of the ordinary potential polynomials to \eqref{eq5}, we obtain
\begin{align}
c_n & = \frac{1}{{\alpha a_0^{\left( {n + \beta } \right)/\alpha } }}\sum\limits_{k = 0}^n {b_{n - k} \mathsf{A}_{ - \left( {n + \beta } \right)/\alpha ,k} \left( {\frac{{a_1 }}{{a_0 }},\frac{{a_2 }}{{a_0 }}, \ldots ,\frac{{a_k }}{{a_0 }}} \right)} \label{eq9} \\
& = \frac{1}{{\alpha a_0^{\left( {n + \beta } \right)/\alpha } }}\sum\limits_{k = 0}^n {b_{n - k} \sum\limits_{j = 0}^k { \binom {- \frac{{n + \beta }}{\alpha }}{j}\frac{1}{{a_0^j }}\mathsf{B}_{k,j} \left( {a_1 ,a_2 , \ldots ,a_{k - j + 1} } \right)} }. \nonumber
\end{align}
This is essentially Wojdylo's formula. To make it more convenient to use, we apply the identity
\[
\binom{ - \rho }{j}= \left( { - 1} \right)^j \binom{j + \rho  - 1}{j}  = \left( { - 1} \right)^j \frac{{\varGamma \left( {j + \rho } \right)}}{{j!\varGamma \left( \rho  \right)}},
\]
which yields
\begin{equation}\label{eq13}
c_n = \frac{1}{{\alpha \varGamma \left( {\frac{{n + \beta }}{\alpha }} \right)}}\sum\limits_{k = 0}^n {b_{n - k} \sum\limits_{j = 0}^k {\frac{{\left( { - 1} \right)^j }}{{a_0^{\left( {n + \beta } \right)/\alpha  + j} }}\frac{{\mathsf{B}_{k,j} \left( {a_1 ,a_2 , \ldots ,a_{k - j + 1} } \right)}}{{j!}}\varGamma \left( {\frac{{n + \beta }}{\alpha } + j} \right)} }.
\end{equation}
This is an even more compact form of Wojdylo's formula. In the common case, when $g \equiv 1$ (and therefore $\beta = 1$), formula \eqref{eq13} is reduced to
\begin{equation}\label{eq17}
c_n  = \frac{1}{{\alpha \varGamma \left( {\frac{{n + 1}}{\alpha }} \right)}}\sum\limits_{k = 0}^n {\frac{{\left( { - 1} \right)^k }}{{a_0^{\left( {n + 1} \right)/\alpha  + k} }}\frac{{\mathsf{B}_{n,k} \left( {a_1 ,a_2 , \ldots ,a_{n - k + 1} } \right)}}{{k!}}\varGamma \left( {\frac{{n + 1}}{\alpha } + k} \right)} .
\end{equation}
Wojdylo's formula \eqref{eq13}, together with the recurrence \eqref{eq32} for the partial ordinary Bell polynomials, provides a systematic way to calculate the coefficients $c_n$.

Comtet \cite[p. 142]{Comtet} showed that an ordinary potential polynomial can be written explicitly in terms of the values of the polynomial at non-negative integers. His formula takes the form
\[
\mathsf{A}_{ - z,k} \left( {\frac{{a_1 }}{{a_0 }},\frac{{a_2 }}{{a_0 }}, \ldots ,\frac{{a_k }}{{a_0 }}} \right) = \frac{{\varGamma \left( {z + k + 1} \right)}}{{k!\varGamma \left( z \right)}}\sum\limits_{j = 0}^k {\frac{\left( { - 1} \right)^j}{z + j} \binom{k}{j}\mathsf{A}_{j,k} \left( {\frac{{a_1 }}{{a_0 }},\frac{{a_2 }}{{a_0 }}, \ldots ,\frac{{a_k }}{{a_0 }}} \right)} .
\]
Applying this in \eqref{eq9}, we obtain the following result.
\begin{corollary}\label{maintheorem}
The coefficients $c_n$ appearing in \eqref{eq2} are given explicitly by
\begin{gather}\label{eq8}
\begin{split}
c_n & = \frac{1}{{\alpha \varGamma \left( {\frac{{n + \beta }}{\alpha }} \right)}}\sum\limits_{k = 0}^n { \frac{\varGamma \left( {\frac{{n + \beta }}{\alpha } + k + 1} \right)b_{n - k}}{k! a_0^{\left( {n + \beta } \right)/\alpha }}\sum\limits_{j = 0}^k {\frac{\left( { - 1} \right)^j}{\frac{{n + \beta }}{\alpha } + j} \binom{k}{j}\mathsf{A}_{j,k}}}\\
& = \frac{1}{{\alpha a_0^{\left( {n + \beta } \right)/\alpha } }}\sum\limits_{k = 0}^n {\binom{ - \frac{{n + \beta }}{\alpha }}{k}  b_{n - k} \sum\limits_{j = 0}^k {\left( { - 1} \right)^{k + j} \frac{{n + \beta  + \alpha k}}{{n + \beta  + \alpha j}} \binom{k}{j}\mathsf{A}_{j,k}} },
\end{split}
\end{gather}
where $\mathsf{A}_{j,k} = \mathsf{A}_{j,k} \left( {a_1 a_0^{-1},a_2 a_0^{-1},\ldots,a_k a_0^{-1}} \right)$.
\end{corollary}
In the common case when $g \equiv 1$, we have $\beta = 1$ and formula \eqref{eq8} is reduced to
\begin{gather}\label{eq14}
\begin{split}
c_n & = \frac{{\varGamma \left( {\frac{{n + 1 }}{\alpha } + n + 1} \right)}}{{\alpha a_0^{\left( {n + 1 } \right)/\alpha } n!\varGamma \left( {\frac{{n + 1 }}{\alpha }} \right)}}\sum\limits_{k = 0}^n {\frac{{\left( { - 1} \right)^k }}{{\frac{{n + 1 }}{\alpha } + k}}\binom{n}{k}\mathsf{A}_{k,n} }\\ 
& = \frac{1}{{\alpha a_0^{\left( {n + 1 } \right)/\alpha } }} \binom{- \frac{{n + 1 }}{\alpha }}{n}\sum\limits_{k = 0}^n {\left( { - 1} \right)^{n + k} \frac{{n + 1  + \alpha n}}{{n + 1  + \alpha k}}\binom{n}{k}\mathsf{A}_{k,n} } .
\end{split}
\end{gather}
From the recurrence relation \eqref{eq33} for the ordinary potential polynomials it follows that
\[
\mathsf{A}_{j,k}  = \sum\limits_{i = 0}^k {\frac{{a_i }}{{a_0 }}\mathsf{A}_{j - 1,k - i} } .
\]
This, together with \eqref{eq8}, provides a simple and efficient way to compute the coefficients $c_n$.

In the next section we shall use the following connection formula between the partial ordinary Bell polynomials and the ordinary potential polynomials (cf. Comtet \cite[p. 156, Exercise 3]{Comtet}).
\begin{equation}\label{eq16}
\mathsf{B}_{k,j}  =  \left( { - 1} \right)^j a_0^j\sum\limits_{i = 0}^j {\left( { - 1} \right)^i \binom{j}{i} \mathsf{A}_{i,k} } .
\end{equation}

\section{Examples}\label{section3}

In this section, we give two illustrative examples to demonstrate the application of Corollary \ref{maintheorem} as well as to compare the results with that given by Wojdylo's formula. In the first example, we derive explicit expressions for the so-called Stirling coefficients appearing in the asymptotic expansion of the gamma function. In the second example, we obtain an explicit formula for certain polynomials related to the coefficients in the uniform asymptotic expansion of the incomplete gamma function.

\subsection*{Example 1} The gamma function can be defined by the following integral
\begin{equation}\label{eq19}
\varGamma \left( {\lambda  + 1} \right) = \int_0^{ + \infty } {\mathrm{e}^{ - t} t^\lambda  \mathrm{d}t}, \quad \lambda >0.
\end{equation}
If we put $t = \lambda\left(1 + x\right)$, we obtain
\[
\varGamma \left( {\lambda  + 1} \right) = \lambda ^{\lambda  + 1} \mathrm{e}^{ - \lambda } \int_{ - 1}^{ + \infty } {\mathrm{e}^{ - \lambda \left( {x - \log \left( {1 + x} \right)} \right)} \mathrm{d}x} ,
\]
and hence, using the identity $\varGamma \left( {\lambda  + 1} \right) = \lambda \varGamma \left( \lambda \right)$,
\begin{equation}\label{eq10}
\frac{{\varGamma \left( \lambda  \right)}}{{\lambda ^\lambda  \mathrm{e}^{ - \lambda } }} = \int_0^{ + \infty } {\mathrm{e}^{ - \lambda \left( {x - \log \left( {1 + x} \right)} \right)} \mathrm{d}x}  + \int_0^1 {\mathrm{e}^{ - \lambda \left( { - x - \log \left( {1 - x} \right)} \right)} \mathrm{d}x}
\end{equation}
follows. Consider the first integral. Let $f\left(x\right) =  x - \log (1 + x)$, $x \geq 0$. Then $f$ has a global minimum at $x=0$ and the expansion
\[
f\left( x \right) = \sum\limits_{k = 0}^\infty  {\left( { - 1} \right)^k \frac{{x^{k + 2} }}{{k + 2}}} 
\]
holds as $x \to 0^+$. We can apply Theorem \ref{theorem1} with $g \equiv 1$, $\alpha = 2$, $\beta = 1$, $b_0 = 1$, $b_k = 0$ for $k>0$, and $a_k = \left( { - 1} \right)^k /\left(k+2\right)$. The result is
\begin{equation}\label{eq28}
\int_0^{ + \infty } {\mathrm{e}^{ - \lambda \left( {x - \log \left( {1 + x} \right)} \right)} \mathrm{d}x} \sim \sum\limits_{n = 0}^\infty  {\varGamma \left( {\frac{{n + 1}}{2}} \right)\frac{{c_n }}{{\lambda ^{\left( {n + 1} \right)/2} }}} ,
\end{equation}
where, by Perron's formula,
\begin{equation}\label{eq29}
c_n  = \frac{1}{{2n!}}\left[ {\frac{{\mathrm{d}^n }}{{\mathrm{d}x^n }}\left( {\frac{{x^2 }}{{x - \log \left( {1 + x} \right)}}} \right)^{\left( {n + 1} \right)/2} } \right]_{x = 0} .
\end{equation}
Similarly, one finds that
\[
\int_0^1 {\mathrm{e}^{ - \lambda \left( {- x - \log \left( {1 - x} \right)} \right)} \mathrm{d}x} \sim \sum\limits_{n = 0}^\infty  {\varGamma \left( {\frac{{n + 1}}{2}} \right)\frac{{\tilde{c}_n }}{{\lambda ^{\left( {n + 1} \right)/2} }}} ,
\]
where
\begin{multline*}
\tilde{c}_n  = \frac{1}{{2n!}}\left[ {\frac{{\mathrm{d}^n }}{{\mathrm{d}x^n }}\left( {\frac{{x^2 }}{{- x - \log \left( {1 - x} \right)}}} \right)^{\left( {n + 1} \right)/2} } \right]_{x = 0}\\ = \frac{\left(-1\right)^n}{{2n!}}\left[ {\frac{{\mathrm{d}^n }}{{\mathrm{d}x^n }}\left( {\frac{{x^2 }}{{x - \log \left( {1 + x} \right)}}} \right)^{\left( {n + 1} \right)/2} } \right]_{x = 0} = \left(-1\right)^n c_n.
\end{multline*}
By substituting these asymptotic series into \eqref{eq10} and performing a simple rearrangement, we deduce
\begin{equation}\label{eq20}
\varGamma \left( \lambda  \right) \sim \sqrt {2\pi } \lambda ^{\lambda  - 1/2} \mathrm{e}^{ - \lambda } \sum\limits_{n = 0}^\infty  {\left( { - 1} \right)^n \frac{{\gamma _n }}{{\lambda ^n }}} 
\end{equation}
as $\lambda \to +\infty$. Here
\begin{equation}\label{eq18}
\gamma _n  = \left( { - 1} \right)^n \sqrt {\frac{2}{\pi }} \varGamma \left( {n + \frac{1}{2}} \right)c_{2n}  = \frac{{\left( { - 1} \right)^n }}{{2^n n!}}\left[ {\frac{{\mathrm{d}^{2n} }}{{\mathrm{d}x^{2n} }}\left( {\frac{1}{2}\frac{{x^2 }}{{x - \log \left( {1 + x} \right)}}} \right)^{n + 1/2} } \right]_{x = 0} 
\end{equation}
are the so-called Stirling coefficients \cite[p. 26]{Paris2}. This representation is well known (see, e.g., \cite{Brassesco}\cite[p. 111]{Wong}). The first few are $\gamma_0 = 1$ and
\[
\gamma _1  =  - \frac{1}{12},\quad \gamma _2  = \frac{1}{288},\quad \gamma _3  = \frac{139}{51840},\quad \gamma _4  =  - \frac{571}{2488320}.
\]
In this case, we are able to obtain explicit formulas for both $\mathsf{B}_{n,k}$ and $\mathsf{A}_{k,n}$. Due to the simpler generating function, we use our new formula \eqref{eq14} and compute the ordinary potential polynomials. Here the generating function takes the form
\[
\left( {2\frac{{x - \log \left( {1 + x} \right)}}{{x^2}}} \right)^k  = \sum\limits_{n = 0}^\infty  {\mathsf{A}_{k,n} x^n } ,
\]
where $\mathsf{A}_{k,n} = \mathsf{A}_{k,n} \left( {a_1 a_0^{-1},a_2 a_0^{-1},\ldots,a_n a_0^{-1}} \right)$ and $a_n = \left( { - 1} \right)^n /\left(n+2\right)$. Using the generating function of the (signless) Stirling numbers of the first kind yields
\begin{align*}
\left( {x - \log \left( {1 + x} \right)} \right)^k & = \sum\limits_{j = 0}^k {\binom{k}{j}x^{k - j} \left( { - \log \left( {1 + x} \right)} \right)^j } \\
& = \sum\limits_{j = 0}^k {\binom{k}{j}x^{k - j} \sum\limits_{n = 0}^\infty  {\left( { - 1} \right)^n j!s\left( {n,j} \right)\frac{{x^n }}{{n!}}} } \\
& = \sum\limits_{n = 0}^\infty  {\left( {\sum\limits_{j = 0}^k {\left( { - 1} \right)^{n - k + j}\binom{k}{j} j!\frac{{s\left( {n - k + j,j} \right)}}{{\left( {n - k + j} \right)!}}} } \right)x^n } ,
\end{align*}
which gives
\begin{equation}\label{eq15}
\mathsf{A}_{k,n} = 2^k \sum\limits_{j = 0}^k {\left( { - 1} \right)^{n + k + j}\binom{k}{j} j!\frac{{s\left( {n + k + j,j} \right)}}{{\left( {n + k + j} \right)!}}} .
\end{equation}
Substituting this into \eqref{eq14} produces
\begin{align*}
c_{2n} &  = \frac{{2^{n-1/2}\varGamma \left( {3n + \frac{3}{2}} \right)}}{{ \left( {2n} \right)!\varGamma \left( {n + \frac{1}{2}} \right)}}\sum\limits_{k = 0}^{2n} {\frac{{2^k }}{{n + \frac{1}{2} + k}}\binom{2n}{k} \sum\limits_{j = 0}^k {\left( { - 1} \right)^{j}\binom{k}{j}j!\frac{{s\left( {2n + k + j,j} \right)}}{{\left( {2n + k + j} \right)!}}}} \\
& = \frac{1}{{\varGamma \left( {n + \frac{1}{2}} \right)}}\sum\limits_{k = 0}^{2n} {\frac{{2^{n + k + 1/2} \varGamma \left( {3n + \frac{3}{2}} \right)}}{{\left( {2n + 2k + 1} \right)\left( {2n - k} \right)!}}\sum\limits_{j = 0}^k {\frac{{\left( { - 1} \right)^j s\left( {2n + k + j,j} \right)}}{{\left( {k - j} \right)!\left( {2n + k + j} \right)!}}} } .
\end{align*}
Finally, by \eqref{eq18},
\begin{equation}\label{eq23}
\gamma_n = \sum\limits_{k = 0}^{2n} {\frac{{\left( { - 1} \right)^n 2^{n + k + 1} \varGamma \left( {3n + \frac{3}{2}} \right)}}{{\sqrt \pi  \left( {2n + 2k + 1} \right)\left( {2n - k} \right)!}}\sum\limits_{j = 0}^k {\frac{{\left( { - 1} \right)^j s\left( {2n + k + j,j} \right)}}{{\left( {k - j} \right)!\left( {2n + k + j} \right)!}}} } .
\end{equation}
As far as we know, this representation of the Stirling coefficients is entirely new. Using expression \eqref{eq16} together with \eqref{eq15}, Wojdylo's formula \eqref{eq17} yields the more complicated formula
\begin{equation}\label{eq21}
\gamma _n  = \sum\limits_{k = 0}^{2n} {\frac{{\left( { - 1} \right)^n 2^n \varGamma \left( {n + k + \frac{1}{2}} \right)}}{{\sqrt \pi  }}} \sum\limits_{j = 0}^k {\frac{{2^j }}{{\left( {k - j} \right)!}}\sum\limits_{i = 0}^j {\frac{{\left( { - 1} \right)^i s\left( {2n + j + i,i} \right)}}{{\left( {j - i} \right)!\left( {2n + j + i} \right)!}}} },
\end{equation}
which is the result of L\'{o}pez, Pagola and P\'{e}rez Sinus\'\i a \cite{Lopez1}.

We remark that the substitution $x = \log(t/\lambda)$ in \eqref{eq19} leads to the form
\[
\frac{{\varGamma \left( \lambda  \right)}}{{\lambda ^\lambda  \mathrm{e}^{ - \lambda } }} = \int_0^{ + \infty } {\mathrm{e}^{ - \lambda \left( {\mathrm{e}^x  - x - 1} \right)} \mathrm{d}x}  + \int_0^{ + \infty } {\mathrm{e}^{ - \lambda \left( {\mathrm{e}^{ - x}  + x - 1} \right)}\mathrm{d}x} .
\]
Similar procedures to the one described above give \eqref{eq20} with
\begin{align}
\gamma _n & = \frac{\left( { - 1} \right)^n}{{2^n n!}}\left[ {\frac{{\mathrm{d}^{2n} }}{{\mathrm{d}x^{2n} }}\left( {\frac{1}{2}\frac{{x^2 }}{{\mathrm{e}^x  - x - 1}}} \right)^{n + 1/2} } \right]_{x = 0} \nonumber \\
& = \sum\limits_{k = 0}^{2n} {\frac{{\left( { - 1} \right)^n 2^{n + k + 1} \varGamma \left( {3n + \frac{3}{2}} \right)}}{{\sqrt \pi  \left( {2n + 2k + 1} \right)\left( {2n - k} \right)!}}\sum\limits_{j = 0}^k {\frac{{\left( { - 1} \right)^j S\left( {2n + k + j,j} \right)}}{{\left( {k - j} \right)!\left( {2n + k + j} \right)!}}} } \label{eq24} \\
& = \sum\limits_{k = 0}^{2n} {\frac{{\left( { - 1} \right)^n 2^n \varGamma \left( {n + k + \frac{1}{2}} \right)}}{{\sqrt \pi  }}} \sum\limits_{j = 0}^k {\frac{{2^j }}{{\left( {k - j} \right)!}}\sum\limits_{i = 0}^j {\frac{{\left( { - 1} \right)^i S\left( {2n + j + i,i} \right)}}{{\left( {j - i} \right)!\left( {2n + j + i} \right)!}}} } \label{eq22},
\end{align}
using formula \eqref{eq6}, \eqref{eq14} and \eqref{eq17}, respectively. Here $S(n,k)$ denotes the Stirling numbers of the second kind. The first representation is known (see, e.g., \cite{Brassesco} \cite{De Angelis}). The second formula is a simplified form of the one derived in \cite{Nemes}. To our knowledge, the third one is entirely new. Note the remarkable similarity between expressions \eqref{eq23} and \eqref{eq24}, and expressions \eqref{eq21} and \eqref{eq22}.

\subsection*{Example 2} The incomplete gamma function is traditionally defined by the following integral:
\[
\varGamma \left( {a,x} \right) = \int_x^\infty  {\mathrm{e}^{ - t} t^{a - 1} \mathrm{d}t} ,\quad a > 0,\quad x \ge 0.
\]
For large values of $x$, the function admits the following asymptotic expansion
\[
\varGamma \left( {a,x} \right)x^{1 - a} \mathrm{e}^x  \sim 1 + \frac{{a - 1}}{x} + \frac{{\left( {a - 1} \right)\left( {a - 2} \right)}}{{x^2 }} +  \cdots ,
\]
which is useful only when $a =o(x)$ \cite[p. 179]{NIST}. However, in exponentially improved asymptotics, we need the asymptotic properties of $\varGamma \left( {a,x} \right)$ as $a \to +\infty$ and $x = \lambda a$, where $\lambda \neq 0$ is a constant. (Do not confuse it with the variable of the integral \eqref{eq1}.) Starting from the integral representation
\[
\frac{{\varGamma \left( {a ,x} \right)}}{{\varGamma \left( a  \right)}} = \frac{{\mathrm{e}^{- a \varphi \left( \lambda  \right)} }}{{2\pi \mathrm{i}}}\int_{c - \mathrm{i}\infty }^{c + \mathrm{i}\infty } {\mathrm{e}^{ a \varphi \left( t \right)} \frac{1}{{\lambda  - t}}\mathrm{d}t} ,\quad 0 < c < \lambda ,
\]
\[
\varphi \left( t \right) =  t - \log t - 1,
\]
Temme \cite{Temme1} gave the following uniform asymptotic expansion of the (normalized) incomplete gamma function as $a \to +\infty$
\[
\frac{{\varGamma \left( {a,x} \right)}}{{\varGamma \left( a \right)}} \sim \frac{1}{2}\mathrm{erfc}\left( {\eta \sqrt {\frac{1}{2}a} } \right) + \frac{{\mathrm{e}^{ - \frac{1}{2}a\eta ^2 } }}{{\sqrt {2\pi a} }}\sum\limits_{n = 0}^\infty  {\frac{{C_n \left( \eta  \right)}}{{a^n }}} .
\]
Here $\eta  = \sqrt {2 \varphi \left( \lambda  \right)}$ and $\mathrm{erfc}$ denotes the complementary error function. He gave a recurrence for the coefficients $C_n \left( \eta  \right)$ and showed that they have the general structure
\begin{equation}\label{eq27}
C_n \left( \eta  \right) = \left( { - 1} \right)^n \left( {\frac{{Q_n \left( \mu  \right)}}{{\mu ^{2n + 1} }} - \frac{{2^n \varGamma \left( {n + \frac{1}{2}} \right)}}{{\sqrt \pi  \eta ^{2n + 1} }}} \right),
\end{equation}
where $\mu = \lambda - 1$ and $Q_n$ is a polynomial in $\mu$ of degree $2n$ \cite{Temme2}. The first few of the $Q_n$'s are given by
\begin{align*}
Q_0 \left( \mu  \right) & = 1,\\
Q_1 \left( \mu  \right) & = 1 + \mu  + \frac{1}{{12}}\mu ^2 ,\\
Q_2 \left( \mu  \right) & = 3 + 5\mu  + \frac{{25}}{{12}}\mu ^2  + \frac{1}{{12}}\mu ^3  + \frac{1}{{288}}\mu ^4 .
\end{align*}

As far as we know, no simple explicit formula for these polynomials exists in the literature. Dunster et al. \cite{Dunster} derived the integral representation
\[
C_n \left( \eta  \right) = \frac{{\left( { - 1} \right)^{n+1} \varGamma \left( {n + \frac{1}{2}} \right)}}{{\left( {2\pi } \right)^{3/2} \mathrm{i}}}\oint_{\mathscr{C}} {\frac{{\mathrm{d}z}}{{\left( {z - \lambda } \right)\left( {z - \log z - 1} \right)^{n + 1/2} }}} .
\]
Here, $\mathscr{C}$ is a loop in the $z$ plane enclosing the poles $z = 1$ and $z = \lambda$. The term $\left( {z - \log z - 1} \right)^{n + 1/2}$ is real and positive for $\arg z = 0$ and $z > 1$ and is defined by continuity elsewhere. We suppose that $\lambda \neq 1$ ($\eta \neq 0$) and use the Residue Theorem to obtain
\begin{align*}
& C_n \left( \eta  \right) = \frac{{\left( { - 1} \right)^{n + 1} \varGamma \left( {n + \frac{1}{2}} \right)}}{{\sqrt {2\pi } }}\left( {\mathop {\mathrm{Res}}\limits_{z = 1} \frac{1}{{\left( {z - \lambda } \right)\left( {z - \log z - 1} \right)^{n + 1/2} }} + \frac{1}{{\left( {\lambda  - \log \lambda  - 1} \right)^{n + 1/2} }}} \right)\\
& = \left( { - 1} \right)^n \left( { - \frac{{\varGamma \left( {n + \frac{1}{2}} \right)}}{{\sqrt {2\pi } }}\left[ {\frac{{\mathrm{d}^{2n} }}{{\mathrm{d}z^{2n} }}\left\{ {\frac{1}{{\left( {z - \lambda } \right)}}\left( {\frac{{\left( {z - 1} \right)^2 }}{{z - \log z - 1}}} \right)^{n + 1/2} } \right\}} \right]_{z = 1}  - \frac{{2^n \varGamma \left( {n + \frac{1}{2}} \right)}}{{\sqrt \pi  \eta ^{2n + 1} }}} \right)\\
& = \left( { - 1} \right)^n \left( {\frac{1}{{2^{2n + 1/2} n!}}\left[ {\frac{{\mathrm{d}^{2n} }}{{\mathrm{d}x^{2n} }}\left\{ {\frac{1}{{\left( {\mu  - x} \right)}}\left( {\frac{{x^2 }}{{x - \log \left( {1 + x} \right)}}} \right)^{n + 1/2} } \right\}} \right]_{x = 0}  - \frac{{2^n \varGamma \left( {n + \frac{1}{2}} \right)}}{{\sqrt \pi  \eta ^{2n + 1} }}} \right) .
\end{align*}
Upon comparing this with \eqref{eq27}, we see that 
\[
Q_n \left( \mu  \right) = \frac{{\mu ^{2n + 1} }}{{2^{2n + 1/2} n!}}\left[ {\frac{{\mathrm{d}^{2n} }}{{\mathrm{d}x^{2n} }}\left\{ {\frac{1}{{\left( {\mu  - x} \right)}}\left( {\frac{{x^2 }}{{x - \log \left( {1 + x} \right)}}} \right)^{n + 1/2} } \right\}} \right]_{x = 0} ,
\]
and from Perron's formula \eqref{eq6},
\begin{equation}\label{eq26}
Q_n \left( \mu  \right) = \frac{{\left( {2n} \right)!\mu ^{2n + 1} }}{{2^{2n - 1/2} n!}}c_{2n}  = \sqrt {\frac{2}{\pi }} \varGamma \left( {n + \frac{1}{2}} \right)\mu ^{2n + 1} c_{2n} ,
\end{equation}
where the $c_n$'s are the coefficients appearing in the asymptotic expansion of the integral \eqref{eq1} with $a=0$, $b>0$, $f\left( x \right) = x - \log \left( {1 + x} \right)$ and $g\left( x \right) = \left( {\mu  - x} \right)^{ - 1}$.

Now, using Corollary \ref{maintheorem}, we derive an explicit formula for the polynomials $Q_n$ in terms of the Stirling numbers of the first kind. Near $x=0$, $g\left( x \right) = \sum\nolimits_{k = 0}^\infty  {\mu ^{ - k - 1} x^k }$, therefore $b_k = \mu ^{ - k - 1}$. The corresponding ordinary potential polynomials are given by \eqref{eq15}. Hence, by Corollary \ref{maintheorem},
\[
c_{2n}  = \frac{1}{{\varGamma \left( {n + \frac{1}{2}} \right)}}\sum\limits_{k = 0}^{2n} {\frac{{2^{n + 1/2} \varGamma \left( {n + k + \frac{3}{2}} \right)\mu ^{ - 2n + k - 1} }}{{k!}}\sum\limits_{j = 0}^k {\frac{{\left( { - 1} \right)^j }}{{2n + 2j + 1}}\binom{k}{j}\mathsf{A}_{j,k} } } .
\]
This, together with \eqref{eq26}, yields
\[
Q_n \left( \mu  \right) = \sum\limits_{k = 0}^{2n} {\frac{{2^{n + 1} \varGamma \left( {n + k + \frac{3}{2}} \right)\mu ^k }}{{\sqrt \pi  k!}}\sum\limits_{j = 0}^k {\frac{{\left( { - 1} \right)^j }}{{2n + 2j + 1}}\binom{k}{j}\mathsf{A}_{j,k} } } ,
\]
or $Q_n \left( \mu  \right) = \sum\nolimits_{k = 0}^{2n} {q_k^{\left( n \right)} \mu ^k }$, where, by \eqref{eq15},
\begin{equation}\label{eq31}
q_k^{\left( n \right)}  = \sum\limits_{j = 0}^k {\frac{{ \left( { - 1} \right)^k 2^{n + j + 1} \varGamma \left( {n + k + \frac{3}{2}} \right)}}{{\sqrt \pi  \left( {2n + 2j + 1} \right)\left( {k - j} \right)!}}\sum\limits_{i = 0}^j {\frac{{\left( { - 1} \right)^i s\left( {k + j + i,i} \right)}}{{\left( {j - i} \right)!\left( {k + j + i} \right)!}}} } .
\end{equation}
We remark that $q_{2n}^{\left( n \right)}  = \left( { - 1} \right)^n \gamma _n$, where $\gamma_n$ is the $n$th Stirling coefficient. Naturally, we could use Wojdylo's formula \eqref{eq13} to obtain an alternative representation for the coefficients $c_n$, and thus for the polynomials $Q_n \left( \mu  \right)$ themselves. The corresponding partial ordinary Bell polynomials can be given explicitly using expression \eqref{eq16} together with \eqref{eq15}, the resulting formula is, however, even more elaborate than \eqref{eq31}.

When $\lambda = 1$ ($\eta = 0$), we have $C_n \left( 0 \right) = \lim _{\eta  \to 0} C_n \left( \eta  \right)$ and Temme's expansion takes the form
\[
\frac{{\varGamma \left( {a,a} \right)}}{{\varGamma \left( a \right)}} \sim \frac{1}{2} + \frac{1}{{\sqrt {2\pi a} }}\sum\limits_{n = 0}^\infty  {\frac{{C_n \left( 0 \right)}}{{a^n }}} 
\]
as $a \to +\infty$. However, the value of the limit $\lim _{\eta  \to 0} C_n \left( \eta  \right)$ is not obvious from \eqref{eq27}. To derive an explicit formula for the $C_n \left( 0 \right)$'s, we proceed as follows: for every positive integer $m$
\begin{align*}
\frac{{\varGamma \left( {m,m} \right)}}{{\varGamma \left( m \right)}} = \mathrm{e}^{ - m} \sum\limits_{k = 0}^{m - 1} {\frac{{m^k }}{{k!}}} & =  - \frac{{m^m \mathrm{e}^{ - m} }}{{m!}} + \frac{{m^m \mathrm{e}^{ - m} }}{{m!}}\sum\limits_{k = 0}^m {\binom{m}{k}\frac{{k!}}{{m^k }}}\\
&  =  - \frac{{m^m \mathrm{e}^{ - m} }}{{m!}} + \frac{{m^m \mathrm{e}^{ - m} }}{{m!}}\sum\limits_{k = 0}^m {\binom{m}{k}m\int_0^{ + \infty } {\mathrm{e}^{ - mx} x^k \mathrm{d}x} } \\
& = \frac{{m^m \mathrm{e}^{ - m} }}{{m!}}\left( { - 1 + m\int_0^{ + \infty } {\mathrm{e}^{ - m\left( {x - \log \left( {1 + x} \right)} \right)} \mathrm{d}x} } \right).
\end{align*}
It is well known that
\begin{equation}\label{eq30}
\frac{{m^m \mathrm{e}^{ - m} }}{{m!}} \sim \frac{1}{{\sqrt {2\pi m} }}\sum\limits_{n = 0}^\infty  {\frac{{\gamma _n }}{{m^n }}},
\end{equation}
where $\gamma_n$ denotes the $n$th Stirling coefficient \cite[p. 26]{Paris2}. By \eqref{eq28}, we have
\[
 - 1 + m\int_0^{ + \infty } {\mathrm{e}^{ - m\left( {x - \log \left( {1 + x} \right)} \right)} \mathrm{d}x}  \sim  - 1 + \sqrt m \sum\limits_{n = 0}^\infty  {\varGamma \left( {\frac{{n + 1}}{2}} \right)\frac{{c_n }}{{m^{n/2} }}} ,
\]
where $c_n$ is given by \eqref{eq29}. From these we deduce
\begin{align*}
& \frac{{\varGamma \left( {m,m} \right)}}{{\varGamma \left( m \right)}} \sim  - \frac{1}{{\sqrt {2\pi m} }}\sum\limits_{n = 0}^\infty  {\frac{{\gamma _n }}{{m^n }}}  + \frac{1}{{\sqrt {2\pi } }}\sum\limits_{n = 0}^\infty  {\frac{{\gamma _n }}{{m^n }}} \sum\limits_{n = 0}^\infty  {\varGamma \left( {\frac{{n + 1}}{2}} \right)\frac{{c_n }}{{m^{n/2} }}} \\
& =  - \frac{1}{{\sqrt {2\pi m} }}\sum\limits_{n = 0}^\infty  {\frac{{\gamma _n }}{{m^n }}}  + \frac{1}{{\sqrt {2\pi } }}\sum\limits_{n = 0}^\infty  {\left( {\sum\limits_{k = 0}^{\left\lfloor {n/2} \right\rfloor } {\gamma _k \varGamma \left( {\frac{{n - 2k + 1}}{2}} \right)c_{n - 2k} } } \right)\frac{1}{{m^{n/2} }}} \\
& = \frac{1}{2} - \frac{1}{{\sqrt {2\pi m} }}\sum\limits_{n = 0}^\infty  {\frac{{\gamma _n }}{{m^n }}}  + \frac{1}{{\sqrt {2\pi m} }}\sum\limits_{n = 0}^\infty  {\underbrace {\left( {\sum\limits_{k = 0}^{\left\lfloor {\left( {n + 1} \right)/2} \right\rfloor } {\gamma _k \varGamma \left( {\frac{n}{2} - k + 1} \right)c_{n - 2k + 1} } } \right)}_{\nu _n }\frac{1}{{m^{n/2} }}} .
\end{align*}
The well-known identity $\sum\nolimits_{k = 0}^n {\left( { - 1} \right)^{n - k} \gamma _k \gamma _{n - k} }  = 0$ ($n \geq 1$) \cite[p. 33]{Paris1} gives
\[
\nu _{2n - 1}  = \sum\limits_{k = 0}^n {\gamma _k \varGamma \left( {n - k + \frac{1}{2}} \right)c_{2n - 2k} }  = \sqrt {\frac{\pi }{2}} \sum\limits_{k = 0}^n {\left( { - 1} \right)^{n - k} \gamma _k \gamma _{n - k} }  = 0,
\]
from which we obtain
\[
\frac{{\varGamma \left( {m,m} \right)}}{{\varGamma \left( m \right)}} \sim \frac{1}{2} + \frac{1}{{\sqrt {2\pi m} }}\sum\limits_{n = 0}^\infty  {\left( { - \gamma_n  + \sum\limits_{k = 0}^n {\gamma _k \varGamma \left( {n - k + 1} \right)c_{2n - 2k + 1} } } \right)\frac{1}{{m^n }}} .
\]
From the uniqueness theorem on asymptotic series we finally have
\begin{align*}
C_n \left( 0 \right) & =  - \gamma_n  + \sum\limits_{k = 0}^n {\gamma _k \varGamma \left( {n - k + 1} \right)c_{2n - 2k + 1} } \\ & =  - \frac{1}{3}\gamma _n  + \sum\limits_{k = 0}^{n - 1} {\left(n - k\right)!\gamma _k c_{2n - 2k + 1} } .
\end{align*}
The first few values are given by
\[
C_0 \left( 0 \right)  =  - \frac{1}{3},\quad C_1 \left( 0 \right)  = -\frac{1}{540},\quad C_2 \left( 0 \right)  = \frac{25}{6048},\quad C_3 \left( 0 \right)  =  \frac{101}{155520}.
\]
These are in agreement with those given in \cite[p. 181]{NIST}.

We remark that one can start with the representation
\[
\frac{{\varGamma \left( {m,m} \right)}}{{\varGamma \left( m \right)}} = \frac{{m^m \mathrm{e}^{ - m} }}{{\varGamma \left( m \right)}}\int_0^{ + \infty } {\mathrm{e}^{ - m\left( {x - \log \left( {1 + x} \right)} \right)} \frac{1}{{1 + x}}\mathrm{d}x} ,
\]
which follows from the definition of the incomplete gamma function. The factor before the integral can be expanded into an asymptotic series using \eqref{eq30}. The asymptotic expansion of the integral can be deduced from Erd\'{e}lyi's theorem. In this way, we get a similar representation for the coefficients $C_n \left( 0 \right)$ to the one above. However, in this case $g\left(x\right) = \left(1+x\right)^{-1}$, whereas in the previous one $g \equiv 1$, which leads to a somewhat complicated representation.

\section{Conclusion and future work}

Laplace's method is one of the classical techniques in the theory of asymptotic expansion of real integrals. The coefficients $c_n$ appearing in the resulting asymptotic expansion, arise as the coefficients of a convergent or asymptotic series of a function defined in an implicit form. Traditionally, series inversion and composition have been used to compute these coefficients, which can be extremely complicated. Nevertheless, there are certain formulas of varying degrees of explicitness for $c_n$'s in the literature. One of them is Perron's formula which gives the coefficients in terms of derivatives of an explicit function. The most explicit formula is given by Wojdylo, who rewrote Perron's formula in terms of some combinatorial objects, called the partial ordinary Bell polynomials.

In this paper we have given an alternative way for simplifying Perron's formula. The new representation involves the ordinary potential polynomials.

We have applied the method to the two important examples of the gamma function and the incomplete gamma function. We have obtained new and explicit formulas for the so-called Stirling coefficients appearing in the asymptotic expansion of the gamma function. We have also derived explicit formulas for certain polynomials related to the coefficients in the uniform asymptotic expansion of the incomplete gamma function.

It turns out that formally the same coefficients appear in the asymptotic expansion of contour integrals with a complex parameter. Our new method, hence, can be applied to give explicit formulas for the coefficients arising in the saddle point method \cite[p. 125]{Olver2}\cite[p. 47]{NIST}, the method of steepest descents \cite{Lopez2}\cite[p. 12]{Paris2} and even in the principle of stationary phase \cite{Olver1}\cite[p. 45]{NIST}.

A fruitful direction for research would be to characterize those integrals that can be transformed into the form \eqref{eq1} and have asymptotic expansions where the ordinary potential polynomials, corresponding to the coefficients, are polynomial-time-computable functions. In the case of the gamma function and the incomplete gamma function, Corollary \ref{maintheorem} produced simpler expressions for the coefficients in their asymptotic expansions than that of Wojdylo's. The question arises whether the ordinary potential polynomial representation generically leads to a simpler form, compared to the partial ordinary Bell polynomial representation, or not. If not, an interesting problem would be to characterize the exceptional cases.

\section*{Acknowledgment}

The author would like to thank the two anonymous referees for their thorough, constructive and helpful comments and suggestions on an earlier version of this paper. Especially, he wishes to thank the second referee, for bringing the relevance of Comtet's work to his attention and giving suggestions for future research.

\appendix

\section{Powers of Power Series}\label{appendix}

Let $F\left( x \right) = 1 + \sum\nolimits_{n = 1}^\infty  {f_n x^n }$ be a formal power series. For any nonnegative integer $k$, we define the partial ordinary Bell polynomials $\mathsf{B}_{n,k} \left( {f_1 ,f_2 , \ldots ,f_{n - k + 1} } \right)$ (associated with $F$) by the generating function
\[
\left( {F\left( x \right) - 1} \right)^k  = \left( {\sum\limits_{n = 1}^\infty  {f_n x^n } } \right)^k  = \sum\limits_{n = k}^\infty  {\mathsf{B}_{n,k} \left( {f_1 ,f_2 , \ldots ,f_{n - k + 1} } \right)x^n },
\]
so that $\mathsf{B}_{0,0}  = 1$, $\mathsf{B}_{n,0}  = 0$ ($n\geq1$), $\mathsf{B}_{n,1}  = f_n$ and $\mathsf{B}_{n,2}  = \sum\nolimits_{j = 1}^{n-1} {f_j f_{n - j} }$. From the simple identity, $\left( {F\left( x \right) - 1} \right)^{k + 1}  = \left( {F\left( x \right) - 1} \right)\left( {F\left( x \right) - 1} \right)^k$, we obtain the recurrence relation
\begin{equation}\label{eq32}
\mathsf{B}_{n,k + 1} \left( {f_1 ,f_2 , \ldots ,f_{n - k} } \right) = \sum\limits_{j = 1}^{n - k} {f_j \mathsf{B}_{n - j,k} \left( {f_1 ,f_2 , \ldots ,f_{n - j - k + 1} } \right)} .
\end{equation}
An explicit representation, that follows from the definition, is given by
\[
\mathsf{B}_{n,k} \left( {f_1 ,f_2 , \ldots ,f_{n - k + 1} } \right) = \sum {\frac{{k!}}{{k_1 !k_2 ! \cdots k_{n - k + 1} !}}f_1^{k_1 } f_2^{k_2 }  \cdots f_{n - k + 1}^{k_{n - k + 1} } },
\]
where the sum runs over all sequences $k_1,k_2, \ldots , k_n$ of non-negative integers such that $k_1  + 2k_2  +  \cdots  + \left(n-k+1\right)k_{n-k+1}  = n$ and $k_1  + k_2  +  \cdots  + k_{n-k+1}  = k$.

For any complex number $\rho$, we define the ordinary potential polynomials $\mathsf{A}_{\rho ,n} \left( {f_1 ,f_2 , \ldots ,f_n } \right)$ (associated to $F$) by the generating function
\[
\left( {F\left( x \right)} \right)^\rho   = \left( {1 + \sum\limits_{n = 1}^\infty  {f_n x^n } } \right)^\rho   = \sum\limits_{n = 0}^\infty  {\mathsf{A}_{\rho ,n} \left( {f_1 ,f_2 , \ldots ,f_n } \right)x^n },
\]
hence,
\[
\mathsf{A}_{\rho ,n} \left( {f_1 ,f_2 , \ldots ,f_n } \right) = \sum\limits_{k = 0}^n {\binom{\rho}{k} \mathsf{B}_{n,k} \left( {f_1 ,f_2 , \ldots ,f_{n - k + 1} } \right)}.
\]
The first few are $\mathsf{A}_{\rho ,0}  = 1$, $\mathsf{A}_{\rho ,1}  = \rho f_1$, $\mathsf{A}_{\rho ,2}  = \rho f_2  + \binom{\rho}{2}f_1^2$, and in general
\[
\mathsf{A}_{\rho ,n} \left( {f_1 ,f_2 , \ldots ,f_n } \right) = \sum {\binom{\rho}{k}\frac{{k!}}{{k_1 !k_2 ! \cdots k_n !}}f_1^{k_1 } f_2^{k_2 }  \cdots f_n^{k_n } } ,
\]
where the sum extends over all sequences $k_1,k_2, \ldots , k_n$ of non-negative integers such that $k_1  + 2k_2  +  \cdots  + nk_n  = n$ and $k_1  + k_2  +  \cdots  + k_n  = k$. Since $\left( {F\left( x \right)} \right)^\rho   = \left( {F\left( x \right)} \right)^{\rho  - 1} F\left( x \right)$, we have the recurrence
\begin{equation}\label{eq33}
\mathsf{A}_{\rho ,n} \left( {f_1 ,f_2 , \ldots ,f_n } \right) = \mathsf{A}_{\rho  - 1,n} \left( {f_1 ,f_2 , \ldots ,f_n } \right) + \sum\limits_{k = 1}^n {f_k \mathsf{A}_{\rho  - 1,n - k} \left( {f_1 ,f_2 , \ldots ,f_{n - k} } \right)} .
\end{equation}
In general, if $G\left( x \right) = \sum\nolimits_{n = 0}^\infty  {g_n x^n }$ is a formal power series, then by definition, we have
\[
G\left( {y\left( {F\left( x \right) - 1} \right)} \right) = \sum\limits_{n = 0}^\infty  {\left( {\sum\limits_{k = 0}^n {g_k \mathsf{B}_{n,k} \left( {f_1 ,f_2 , \ldots ,f_{n - k + 1} } \right)y^k } } \right)x^n } .
\]
Specially,
\[
\exp \left( {y\sum\limits_{n = 1}^\infty  {f_n x^n } } \right) = \sum\limits_{n = 0}^\infty  {\left( {\sum\limits_{k = 0}^n {\frac{{\mathsf{B}_{n,k} \left( {f_1 ,f_2 , \ldots ,f_{n - k + 1} } \right)}}{{k!}}y^k } } \right)x^n } .
\]
For more details see, e.g., Comtet's book \cite[pp. 133--153]{Comtet}.

\end{document}